\documentclass[3p]{elsarticle}

\usepackage{amsthm}
\usepackage{amssymb}
\usepackage{amsmath}
\usepackage{amsfonts}
\usepackage{epsfig}
\usepackage[mathscr]{eucal}
\usepackage[ruled,vlined]{algorithm2e}
\usepackage{algorithmic}
\usepackage{setspace,color}
\usepackage{tikz}
\usepackage{hyperref}
\usepackage[capitalize]{cleveref}
\usepackage{numcompress}
\usepackage[caption=false]{subfig}
\usepackage{multirow}

\usepackage[numbers]{natbib}
\usepackage{stackrel}

\usetikzlibrary{snakes}
\usetikzlibrary{arrows,shapes}

\def\R{{\mathbb R}}
\def\C{{\mathbb C}}

\def\N{{\mathbb N}}

\def\Z{{\mathbb Z}}

\def\PP{\mathbb{P}}

\def\su{\mathrm{supp}}

\def\TV{\mathrm{TV}}

\def\Om{\Omega}

\def\ep{\epsilon}

\def\f{\frac}
\def\p{\partial}

\def\na{\nabla}
\def\la{\langle}
\def\ra{\rangle}

\def\rd{{\mathrm d}}
\def\e{{\boldsymbol e}}

\def\bsf{{\boldsymbol f}}

\def\bsg{{\boldsymbol g}}

\def\bi{{\mathbf i}}

\def\bk{{\boldsymbol k}}
\def\bsl{{\boldsymbol l}}

\def\bq{\boldsymbol{q}}

\def\bu{\boldsymbol{u}}

\def\x{\boldsymbol{x}}
\def\y{{\boldsymbol y}}
\def\z{{\boldsymbol z}}

\def\I{{\mathbf I}}

\def\bfQ{{\mathbf Q}}

\def\bS{{\mathbf S}}

\def\mE{{\mathcal E}}

\def\mN{{\mathcal N}}

\def\mR{{\mathcal R}}

\def\bmR{{\boldsymbol \mR}}

\def\msM{{\mathscr M}}
\def\msN{{\mathscr N}}

\def\msR{{\mathscr R}}

\def\msT{{\mathscr T}}

\def\eep{\boldsymbol{\eps}}

\def\oom{\boldsymbol{\omega}}
\def\ssi{\boldsymbol{\sigma}}

\def\bi{\begin{itemize}} \def\ei{\end{itemize}}
\def\be{\begin{eqnarray*}}
\def\ee{\end{eqnarray*}}

\def\0{{\mathbf 0}}

\newcommand{\beq}{\begin{equation}}
\newcommand{\eeq}{\end{equation}}

\def\La{{\boldsymbol\Lambda}}

\def\wt{\widetilde}
\def\wh{\widehat}

\def\Na{\boldsymbol \nabla}

\def\bOm{\boldsymbol{\Omega}}

\newcommand{\eps}{\varepsilon}
\def\la{\langle}
\def\ra{\rangle}

\def\XXint#1#2#3{{\setbox0=\hbox{$#1{#2#3}{\int}$ }
\vcenter{\hbox{$#2#3$ }}\kern-.55\wd0}}

\newtheorem{thm}{Theorem}[section]
\newtheorem{proposition}[thm]{Proposition}
\newtheorem{lemma}[thm]{Lemma}
\newtheorem{theorem}[thm]{Theorem}

\newdefinition{definition}[thm]{Definition}
\newdefinition{rmk}[thm]{Remark}
\newdefinition{notation}[thm]{Notation}
\newdefinition{example}[thm]{Example}

\newproof{pf}{Proof}
\newproof{poth4}{Proof of \cref{Th4}}

\newcommand{\argmin}{\operatornamewithlimits{argmin}}

\numberwithin{equation}{section}

\crefformat{equation}{\textup{#2(#1)#3}}
\crefrangeformat{equation}{\textup{#3(#1)#4--#5(#2)#6}}
\crefmultiformat{equation}{\textup{#2(#1)#3}}{ and \textup{#2(#1)#3}}
{, \textup{#2(#1)#3}}{, and \textup{#2(#1)#3}}
\crefrangemultiformat{equation}{\textup{#3(#1)#4--#5(#2)#6}}%
{ and \textup{#3(#1)#4--#5(#2)#6}}{, \textup{#3(#1)#4--#5(#2)#6}}{, and \textup{#3(#1)#4--#5(#2)#6}}

\Crefformat{equation}{#2Equation~\textup{(#1)}#3}
\Crefrangeformat{equation}{Equations~\textup{#3(#1)#4--#5(#2)#6}}
\Crefmultiformat{equation}{Equations~\textup{#2(#1)#3}}{ and \textup{#2(#1)#3}}
{, \textup{#2(#1)#3}}{, and \textup{#2(#1)#3}}
\Crefrangemultiformat{equation}{Equations~\textup{#3(#1)#4--#5(#2)#6}}%
{ and \textup{#3(#1)#4--#5(#2)#6}}{, \textup{#3(#1)#4--#5(#2)#6}}{, and \textup{#3(#1)#4--#5(#2)#6}}

\crefdefaultlabelformat{#2\textup{#1}#3}

\title{Approximation Theory of Total Variation Minimization for Data Completion}

\author[HKUST]{Jian-Feng Cai\fnref{JFC}}
\ead{jfcai@ust.hk}
\author[TJU]{Jae Kyu Choi\fnref{JKC}\corref{cor}}
\ead{jaycjk@tongji.edu.cn}
\author[FUDAN]{Ke Wei}
\ead{kewei@fudan.edu.cn}

\cortext[cor]{Corresponding author}

\address[HKUST]{Department of Mathematics, The Hong Kong University of Science and Technology, Hong Kong}
\address[TJU]{School of Mathematical Sciences, Tongji University, Shanghai, 200092 China}
\address[FUDAN]{School of Data Science, Fudan University, Shanghai, 200433 China}

\fntext[JFC]{J. F. Cai is supported by the Hong Kong Research Grants Council (HKRGC) GRF 16306821, 16309219, and 16310620.}
\fntext[JKC]{J. K. Choi is supported in part by the National Natural Science Foundation of China (NSFC) Youth Program 11901436, the Shanghai Science and Technology Committee 20JC1413500, and the Fundamental Research Funds for the Central Universities.}
\begin{document}

\begin{abstract} Total variation (TV) minimization is one of the most important techniques in modern signal/image processing, and has wide range of applications. While there are numerous recent works on the restoration guarantee of the TV minimization in the framework of compressed sensing, there are few works on the restoration guarantee of the restoration from partial observations. This paper is to analyze the error of TV based restoration from random entrywise samples. In particular, we estimate the error between the underlying original data and the approximate solution that interpolates (or approximates with an error bound depending on the noise level) the given data that has the minimal TV seminorm among all possible solutions. Finally, we further connect the error estimate for the discrete model to the sparse gradient restoration problem and to the approximation to the underlying function from which the underlying true data comes.
\end{abstract}

\begin{keyword} Total variation \sep missing data restoration \sep error estimation \sep $\ell_1$ minimization \sep uniform law of large numbers \sep covering number
\end{keyword}

\maketitle

\pagestyle{myheadings}
\thispagestyle{plain}
\markboth{Jian-Feng Cai, Jae Kyu Choi, and Ke Wei}{Approximation Theory of TV Minimization for Data Completion}

\section{Introduction}

Total variation (TV) minimization is one of the most important techniques in modern signal/image processing, and it has wide range of applications in denoising, deblurring, and inpainting \cite{B.Adcock2021,A.Chambolle2010,A.Chambolle2016}. In this paper, we consider the restoration from partial observations, which is to restore the underlying data $\bsf$ from a given partial entrywise observation $\bsg$ by solving
\begin{align}\label{Linear_IP_Incomplete}
\bsg[\bk]=\left\{\begin{array}{cl}
\bsf[\bk]+\eep[\bk],~&\bk\in\La,\vspace{0.25em}\\
\text{unknown},~&\bk\in\bOm\setminus\La.
\end{array}\right.
\end{align}
In \cref{Linear_IP_Incomplete}, $\eep$ is measurement error which is assumed to be uncorrelated with $\bsf\big|_{\La}$, $\bOm$ is the domain on the lattice $\Z^d$ ($d\in\N$) where the underlying true data $\bsf$ is defined, and $\La\subseteq\bOm$ is a set of indices where $\bsg$ is reliable. Note that this paper involves both functions and their discrete counterparts. We shall use regular characters to denote functions and use bold-faced characters to denote their discrete analogies. For example, we use $u$ as an element in a function space, while we use $\bu$ to denote its corresponding discretized version (the type of discretization will be made clear later).

The partial observation $\bsg$ can be part of sound, images, time-varying measurement values and sensor data, and the goal is to fill-in the missing region $\bOm\setminus\La$ from the measurement $\bsg$ \cite{J.F.Cai2011}. Related tasks include e.g. \cite{M.Bertalmio2003,A.Bugeau2010,J.F.Cai2008b,T.F.Chan2002,M.Elad2005} for image inpainting, \cite{J.F.Cai2010b,E.J.Candes2010,E.J.Candes2009} for matrix completion, \cite{F.Cucker2002,Vapnik1998} for regression in machine learning, \cite{M.J.Johnson2009} for surface reconstruction in computer graphics, and \cite{E.J.Candes2006,R.H.Chan2005} for miscellaneous applications. Nevertheless, we forgo to give a detailed survey on these various applications and the interested reader should consult the references mentioned above for details. Instead, assuming that
\begin{align}\label{Varianceofg}
\mathrm{Var}(\eep)\leq\eta^2\leq\mathrm{Var}(\bsg):=\f{1}{|\La|}\sum_{\bk\in\La}\left|\bsg[\bk]-\f{1}{|\La|}\sum_{\bsl\in\La}\bsg[\bsl]\right|^2<\infty,
\end{align}
we focus on the following constrained TV minimization problem
\begin{align}\label{TVMinimization}
\min_{\bu}~\left\|\Na\bu\right\|_1~~~\text{subject to}~~\f{1}{|\La|}\sum_{\bk\in\La}\left|\bu[\bk]-\bsg[\bk]\right|^2\leq\eta^2.
\end{align}
In \cref{TVMinimization}, the TV seminorm $\left\|\Na\bu\right\|_1$ is defined as
\begin{align}\label{AnisoTV}
\left\|\Na\bu\right\|_1=\sum_{j=1}^d\sum_{\bk,\bk+\e_j\in\bOm}\left|\bu[\bk+\e_j]-\bu[\bk]\right|,
\end{align}
where $\left\{\e_1,\cdots,\e_d\right\}$ denotes the standard basis of $\R^d$.

In this paper, we will consider the following random setting. For
\begin{align}\label{OmegaSetting}
\bOm=\left\{0,1,\cdots,N-1\right\}^d
\end{align}
with $N\in\N$, let
\begin{align}\label{LambdaSetting}
\La\subseteq\bOm,~~|\La|=m,~~\La~\text{is uniformly random chosen from all $m$-subsets of}~\bOm.
\end{align}
To simplify the discussion, we only consider the regular $d$ dimensional hypercubic grid $\bOm$. Note, however, that it is not difficult to generalize our analysis to the arbitrary $d$ dimensional hyperrectangular grid. Then we denote by $\rho:=m/|\bOm|$ the density of the pixels available. Notice that we assume that the measurement $\bsg$ and the error $\eep$ are given and fixed, even though the error $\eep$ can be viewed as a particular realization of some random variables, e.g. the i.i.d. Gaussian noise. Hence, the only random variable in our setting is the data set $\La$ which is uniformly drawn from all $m$-subsets of $\bOm$. Such a missing data restoration from randomly sampled pixels frequently occur when part of pixels is randomly missing due to e.g. the unreliable communication channel \cite{T.F.Chan2006a}, and/or the corruption by a salt-and-pepper noise \cite{J.F.Cai2009,R.H.Chan2005}.

The focus of this paper is to study the approximation property of the above TV minimization problem \cref{TVMinimization}. To do this, we assume that the underlying true data $\bsf$ satisfies $\left\|\Na\bsf\right\|_1<\infty$ and $\bsf[\bk]\in[0,M]$, $\bk\in\bOm$, for some constant $M\geq1$. The first condition enforces some regularity on $\bsf$ and the second condition enforces the boundedness of each pixel value. In addition, since the measurement $\bsg$ can be ``clipped'' due to e.g. the hardware limitation \cite{J.A.Tropp2010}, we further assume that $\bsg$ satisfies $\bsg[\bk]\in[0,M]$ for $\bk\in\La$ as well. Let $\bu^{\La}$ be a solution to \cref{TVMinimization}. Notice that, as a discrete analogy of \cite{T.F.Chan2002}, it is not hard to verify that \cref{TVMinimization} admits at least one solution with the minimal total variation seminorm subject to the constraint. In addition, it is obvious that, if $\rho=1$ (i.e. $\La=\bOm$) and $\eta=0$, $\bsf$ is the unique solution. Hence, we are interested in analyzing what happens when $\rho<1$. More precisely, we will show that, under some mild assumptions, the error between $\bu^{\La}$ and $\bsf$ satisfies
\begin{align}\label{Goal}
\f{1}{|\bOm|}\left\|\bu^{\La}-\bsf\right\|_{\ell_2(\bOm)}^2\leq C\rho^{-1/2}|\bOm|^{-\beta}\left(\log_2|\bOm|\right)^{3/2}+\f{16}{3}\eta^2
\end{align}
with probability at least $1-|\bOm|^{-1}$. In \cref{Goal}, $\beta>0$ is a constant related to the regularity of $\bsf$, and $C>0$ is a constant independent of $\rho$, $|\bOm|$ (or equivalently $N$), and $\eta$. Roughly speaking, as long as $\bOm$ is sufficiently large, there exists a pretty good chance to restore a data close to the underlying original one $\bsf$ by solving \cref{TVMinimization} from $m$ partial observation $\bsg$.

In the literature, the restoration guarantees of TV minimization has been extensively studied in the various signal/image restoration tasks. For instance, the authors in \cite{E.J.Candes2006} asserted the restoration guarantee of $\bu\in\C^N$ with $s$-sparse gradient from $m$ noiseless Fourier samples drawn uniformly and randomly. Briefly speaking, as long as the number of samples $m$ satisfies $m\gtrsim s\ln N$, $\bu$ can be exactly restored with high probability by solving the TV minimization with the equality constraint on the Fourier samples. Later, the authors in \cite{D.Needell2013,D.Needell2013a} further extended the result to the $d\geq2$ dimensional TV image restoration problems. More precisely, the authors present the connection between the TV seminorm and the compressibility under the Haar wavelet (e.g. \cite{Daubechies1992,Mallat2008}) to modify a restricted isometry property (RIP) of a random sensing matrix under the Haar wavelet representation. Unlike the works in \cite{D.Needell2013,D.Needell2013a} with the RIP condition, the authors in \cite{J.F.Cai2015} present a restoration guarantee based on the nullspace condition of a Gaussian sensing matrix. Briefly speaking, the authors use the so-called ``Escape through the Mesh'' theorem (e.g. \cite{V.Chandrasekaran2012}) to estimate the Gaussian width (e.g. \cite{M.Rudelson2008}) of a cone
specified by the null space property (NSP) condition. Finally, apart from the aforementioned works on the random sensing matrices, the restoration guarantees for Fourier samples have also been extensively studied in \cite{B.Adcock2021,F.Krahmer2014,Poon2015} under various sampling strategies. While these aforementioned restoration guarantees are mainly in the context of compressed sensing, to the best of our knowledge, there are few works on the error analysis of the missing data restoration from partial measurements.

The analysis in this paper is different from the aforementioned previous works. First of all, our analysis does not require explicit sparse gradient. Notice that the ``sensing matrix'' of \cref{Linear_IP_Incomplete} does not satisfy the RIP condition or the NSP condition for the sparse restoration. Instead, we assume the bounded TV seminorm to impose some mild regularity condition of the underlying image, rather than the sparse gradient, so our analysis is not restricted to the data with sparse gradient. In fact, our analysis mostly follows the direction of \cite{J.F.Cai2011} on the approximation property of frame based missing data restoration. More precisely, we also use the combination of the uniform law of large numbers, which is standard in classical empirical processes and statistical learning theory, and an estimation for its involved covering number of a hypothesis space of the solution. Nevertheless, we further mention that our error analysis uses a new estimation for covering number. In other words, even though we also use the special structure of the set and the max-flow min-cut theorem in graph theory to estimate covering number, we improve the estimate in \cite[Theorem 2.4]{J.F.Cai2011} by relaxing the constraint of the radius to an arbitrarily small one. As a consequence, our analysis is not limited to \cref{TVMinimization}. In fact, with the aid of our new estimate for covering number, it is not difficult to extend our analysis to the analysis of tight frame based missing data restoration in \cite{J.F.Cai2011} with a slight modification.

The rest of this paper is organized as follows. In \cref{ErrorAnalysis}, we give our main results of approximation analysis for TV minimization \cref{TVMinimization}. In \cref{SparseRestoration,TVInpainting}, we illustrate the applications of our main results. More precisely, we estimate the error of piecewise constant data restoration from random samples in \cref{SparseRestoration}. In \cref{TVInpainting}, we estimate the error of TV image inpainting, and further connect the error analysis for the discrete problem \cref{TVMinimization} to the two dimensional $BV$ function approximation based on the finite element approximation of $BV$ functions \cite{S.Bartels2015}. Finally, all technical proofs are postponed to \cref{TechnicalProofs}, and \cref{Conclusion} concludes this paper with some future directions.

\section{Approximation error of TV minimization}\label{ErrorAnalysis}

In this section, we give the approximation error analysis of the TV minimization model \cref{TVMinimization} with the anisotropic TV seminorm $\left\|\Na\bu\right\|_1$ defined as \cref{AnisoTV}. Specifically, let $\bu^{\La}$ be the solution to \cref{TVMinimization}, and we derive the explicit formulation of \cref{Goal}. Recall that $\bOm=\left\{0,\cdots,N-1\right\}^d$ and $\La$ is a data set which is uniformly randomly chosen from all $m$-subsets of $\bOm$.

To begin with, we need to define an appropriate hypothesis space. Hence, we present \cref{Lemma1} related to the discrete maximum principle (e.g. \cite{Bartels2016}) to characterize the solution $\bu^{\La}$.

\begin{lemma}\label{Lemma1} Assume that $\bsg\in[0,M]^{\La}$ in \cref{Linear_IP_Incomplete} is a nonconstant vector which satisfies \cref{Varianceofg}. For a solution $\bu^{\La}$ to \cref{TVMinimization}, we have
\begin{align}\label{DiscreteMaximum}
0\leq\bu^{\La}[\bk]\leq M~~~~~\text{for all}~~~\bk\in\bOm.
\end{align}
\end{lemma}

\begin{pf} Let $\bu^{\La}$ be a solution to \cref{TVMinimization}. It is obvious that, if we have
\begin{align*}
\mathrm{Var}(\bsg)=\f{1}{|\La|}\sum_{\bk\in\La}\left|\bsg[\bk]-\mE\left(\bsg\right)\right|^2=\eta^2~~~\text{with}~~~\mE\left(\bsg\right)=\f{1}{|\La|}\sum_{\bk\in\La}\bsg[\bk],
\end{align*}
then the constant vector $\bu^{\La}[\bk]=\mE\left(\bsg\right)$ for all $\bk\in\bOm$ is the minimizer of \cref{TVMinimization} with \cref{DiscreteMaximum}. Hence, it suffices to prove that \cref{DiscreteMaximum} holds when
\begin{align*}
\eta^2<\f{1}{|\La|}\sum_{\bk\in\La}\left|\bsg[\bk]-\mE\left(\bsg\right)\right|^2.
\end{align*}
In this case, since we have
\begin{align*}
\f{1}{|\La|}\sum_{\bk\in\La}\left|\bsg[\bk]-c\right|^2\geq\mathrm{Var}(\bsg)>\eta^2,
\end{align*}
for any constant $c\in\R$, \cref{TVMinimization} does not admit a constant vector as a minimizer. For each $\bk\in\bOm$, we define $\wt{\bu}^{\La}$ by
\begin{align}\label{TruncatedMinimizer}
\wt{\bu}^{\La}[\bk]=\left\{\begin{array}{ccl}
\bu^{\La}[\bk]~&\text{if}&\bu^{\La}[\bk]\leq M\vspace{0.25em}\\
M~&\text{if}&\bu^{\La}[\bk]>M.
\end{array}\right.
\end{align}
Obviously, we have $\wt{\bu}^{\La}[\bk]\leq M$. In addition, we can write $\wt{\bu}^{\La}[\bk]=G(\bu^{\La}[\bk])$ for $\bk\in\bOm$, where $G:\R\to\R$ is defined as
\begin{align}\label{FunctionG}
G(x)=\left\{\begin{array}{ccl}
x~&\text{if}&x\leq M_{\bsg}\vspace{0.25em}\\
M~&\text{if}&x>M.
\end{array}\right.
\end{align}
Notice that this $G:\R\to\R$ is Lipschitz continuous with Lipschitz constant $1$:
\begin{align}\label{Lipschitz}
\left|G(x)-G(y)\right|\leq\left|x-y\right|,
\end{align}
and takes $\bsg$ as the fixed point: $G(\bsg[\bk])=\bsg[\bk]$ for $\bk\in\La$.

From \cref{Lipschitz}, we have
\begin{align*}
\left|\wt{\bu}^{\La}[\bk]-\bsg[\bk]\right|=\left|G(\wt{\bu}^{\La}[\bk])-G(\bsg[\bk])\right|\leq\left|\bu^{\La}[\bk]-\bsg[\bk]\right|,
\end{align*}
which means that $\wt{\bu}^{\La}$ satisfies the constraint:
\begin{align*}
\f{1}{|\La|}\sum_{\bk\in\La}\left|\wt{\bu}^{\La}[\bk]-\bsg[\bk]\right|^2\leq\f{1}{|\La|}\sum_{\bk\in\La}\leq\left|\bu^{\La}[\bk]-\bsg[\bk]\right|^2\leq\eta^2.
\end{align*}
In addition, \cref{Lipschitz} also gives us
\begin{align*}
\left\|\Na\wt{\bu}^{\La}\right\|_1&=\sum_{j=1}^d\sum_{\bk,\bk+\e_j\in\bOm}\left|\wt{\bu}^{\La}[\bk+\e_j]-\wt{\bu}^{\La}[\bk]\right|\\
&\leq\sum_{j=1}^d\sum_{\bk,\bk+\e_j\in\bOm}\left|\bu^{\La}[\bk+\e_j]-\bu^{\La}[\bk]\right|=\left\|\Na\bu^{\La}\right\|_1.
\end{align*}
Since $\bu^{\La}$ is a minimizer of \cref{TVMinimization}, $\left\|\Na\bu^{\La}\right\|_1$ should be the minimum discrete total variation subject to the constraint. Hence, it follows that $\left\|\Na\wt{\bu}^{\La}\right\|_1=\left\|\Na\bu^{\La}\right\|_1$, and in particular, by \cref{TruncatedMinimizer}, $\wt{\bu}^{\La}=\bu^{\La}$. If not, then $\bu^{\La}[\bk_0]>M=\wt{\bu}^{\La}[\bk_0]$ for some $\bk_0\in\bOm$, which leads to a contradiction that $\left\|\Na\wt{\bu}^{\La}\right\|_1<\left\|\Na\bu^{\La}\right\|_1$. Hence, it must be that $\bu^{\La}[\bk]=\wt{\bu}^{\La}[\bk]\leq M$ for all $\bk\in\bOm$. For the lower bound, we use the same argument by considering $\bsg\mapsto-\bsg$, i.e. noting that $-\bu^{\La}$ is a minimizer of \cref{TVMinimization}. This completes the proof.\qquad$\square$
\end{pf}

With the aid of \cref{Lemma1}, we can consider the following set
\begin{align}\label{HypothesisSpace}
\msM=\left\{\bu\in\ell_{\infty}(\bOm):\left\|\Na\bu\right\|_1\leq\left\|\Na\bsf\right\|_1,~\f{1}{|\La|}\sum_{\bk\in\La}\left|\bu[\bk]-\bsg[\bk]\right|^2\leq\eta^2,~\bu\in[0,M]^{\bOm}\right\},
\end{align}
as an involved hypothesis space. In \cref{HypothesisSpace}, $M\geq1$ is a positive constant related to the boundedness of each pixel value, $\eta>0$ is a fixed positive constant related to the bound of measurement error, and $\left\|\Na\bu\right\|_1$ is the anisotropic TV defined as \cref{AnisoTV}. Notice that, given that $\bsg\in[0,M]^{\La}$ is defined as \cref{Linear_IP_Incomplete} with $\bsf\in[0,M]^{\bOm}$, the true solution $\bsf$ lies in this set $\msM$. In addition, since $\bu^{\La}$ is a solution to \cref{TVMinimization}, it follows that $\left\|\Na\bu^{\La}\right\|_1\leq\left\|\Na\bsf\right\|_1$, and $\bu^{\La}\in[0,M]^{\bOm}$ by \cref{Lemma1}. Hence, the set $\msM$ defined as \cref{HypothesisSpace} is the desired space.

Our error analysis relies on the capacity of an involved set. Notice that there are numerous tools including VC dimension \cite{Vapnik1998}, $V_{\gamma}$-dimension and $P_{\gamma}$-dimension \cite{N.Alon1997}, Rademacher complexities \cite{P.L.Bartlett2006,Koltchinskii2001}, and covering number \cite{F.Cucker2002}. In this paper, we choose the covering number to measure the capacity of the hypothesis space, as it is the most convenient tool for metric spaces \cite{J.F.Cai2011}.

\begin{definition} Let $\msM\subseteq\R^{\bOm}$ and $r>0$ be given. The covering number $\mN(\msM,r)$ is defined as
\begin{align*}
\mN\left(\msM,r\right)=\inf\left\{K\in\N:\exists~\bu_1,\cdots,\bu_K\in\msM~~\text{s.t.}~~\msM\subseteq\bigcup_{j=1}^{K}\left\{\bu\in\msM:\left\|\bu-\bu_j\right\|_{\ell_{\infty}(\bOm)}\leq r\right\}\right\}.
\end{align*}
\end{definition}

With this idea of covering number, we present the first relation between the solution $\bu^{\La}$ of \cref{TVMinimization} and the underlying true image $\bsf$. Specifically, we estimate the probability of the event
\begin{align*}
\f{1}{|\bOm|}\left\|\bu^{\La}-\bsf\right\|_{\ell_2(\bOm)}^2\leq\ep+\f{16}{3}\eta^2
\end{align*}
for an arbitrary $\ep>0$ in terms of the covering number of the hypothesis space $\msM$ defined as \cref{HypothesisSpace}. Since the proof is exactly the same as \cite[Theorem 2.3]{J.F.Cai2011}, we omit the proof.

\begin{proposition}\label{Prop3} Let $\msM$ be defined as \cref{HypothesisSpace} and $\bu^{\La}$ be a solution to \cref{TVMinimization}. Then for an arbitrary $\ep>0$, the following inequality
\begin{align}\label{Prob}
\PP\left\{\f{1}{|\bOm|}\left\|\bu^{\La}-\bsf\right\|_{\ell_2(\bOm)}^2\leq\ep+\f{16}{3}\eta^2\right\}\geq1-\mN\left(\msM,\f{\ep}{12M}\right)\exp\left(-\f{3m\ep}{256M^2}\right)
\end{align}
holds. In \cref{Prob}, $m=|\La|$ denotes the number of samples.
\end{proposition}

Hence, our error estimate will be completed if we can bound the covering number in \cref{Prob}. At first glance, for each $\bu\in\msM$, we have $\left\|\bu\right\|_{\ell_{\infty}(\bOm)}\leq M$, so we have the following simple upper bound
\begin{align}\label{CoveringNumberRoughEstimate}
\mN(\msM,r)\leq\left(\f{2M}{r}\right)^{|\bOm|},
\end{align}
as presented in \cite{F.Cucker2002}. However, since the above estimation \cref{CoveringNumberRoughEstimate} is not tight enough to derive an error estimate, we need to find a much tighter upper bound for $\mN(\msM,r)$ by further exploiting the conditions of $\msM$. Notice that in our setting, the bounded TV seminorm is to impose some regularity condition on the image to be restored, whence it would be quite reasonable to get a tighter bound by exploiting the condition $\left\|\Na\bu\right\|_1\leq\left\|\na\bsf\right\|_1$. This leads us to relax $\msM$ in \cref{HypothesisSpace} into
\begin{align}\label{HypothesisSpaceRelax}
\wt{\msM}=\left\{\bu\in\ell_{\infty}(\bOm):\left\|\Na\bu\right\|_1\leq\left\|\Na\bsf\right\|_1,~\left\|\bu\right\|_{\ell_{\infty}(\bOm)}\leq M\right\}.
\end{align}
With this $\wt{\msM}$, we can obtain the desired estimate of the covering numbers in \cref{Th1}. Similar to \cite[Theorem 2.4]{J.F.Cai2011}, our estimate is also based on the quantized total variation minimization (e.g. \cite{A.Chambolle[2021]copyright2021,A.Chambolle2009}) and the max-flow min-cut theorem \cite{Diestel2018}. In this paper, we improve \cite[Theorem 2.4]{J.F.Cai2011} by relaxing the constraint on the radius $r$. Since the proof is long and technical, it is postponed to \ref{PfTh1}.

\begin{theorem}\label{Th1} Let $\msM$ be defined as \cref{HypothesisSpace}. Assume that $\bsf$ satisfies $\left\|\Na\bsf\right\|_1\leq C_{\bsf}\left|\bOm\right|^b$ for some $b\in[0,1]$. Then for $r\geq|\bOm|^{-a}$ with $a\geq1-b$, we have
\begin{align}\label{CoveringNumberNewEstimate}
\ln\mN\left(\msM,r\right)\leq\f{C_{a,b,d}\left|\bOm\right|^b}{r}\log_2\left|\bOm\right|
\end{align}
where $C_{a,b,d}=40(2a+b)MC_{\bsf}\left(d+2C_{\bsf}\right)$.
\end{theorem}

With the aid of \cref{Prop3,Th1}, we can give the explicit form of our main result. Briefly speaking, for a fixed $\rho$, as long as the cardinality of $\bOm$ is sufficiently large, the solution $\bu^{\La}$ to \cref{TVMinimization} gives a good approximation of the original data $\bsf$. Moreover, for a fixed $\bOm$, the error becomes smaller as $\rho\nearrow1$, which coincides with the common sense as larger observations are available for a fixed $\bOm$ with larger $\rho$.

\begin{theorem}\label{Th2} Assume that $\bsg\in[0,M]^{\La}$ in \cref{Linear_IP_Incomplete} is a nonconstant vector which satisfies \cref{Varianceofg}. Let $\bu^{\La}$ be a solution to \cref{TVMinimization}, and let $\left\|\Na\bsf\right\|_1\leq C_{\bsf}\left|\bOm\right|^b$ for some $b\in[0,1)$. Then the following inequality
\begin{align}\label{MainResult}
\f{1}{|\bOm|}\left\|\bu^{\La}-\bsf\right\|_{\ell_2(\bOm)}^2\leq\wt{c}\rho^{-1/2}|\bOm|^{-\f{1-b}{2}}\left(\log_2|\bOm|\right)^{3/2}+\f{16}{3}\eta^2
\end{align}
holds with probability at least $1-|\bOm|^{-1}$. In \cref{MainResult}, $\wt{c}$ is defined as
\begin{align*}
\wt{c}=\f{64}{3}M^2\left(4+3\sqrt{10(2a+b)C_{\bsf}\left(d+2C_{\bsf}\right)}\right),
\end{align*}
for some $a\geq1-b$.
\end{theorem}

\begin{pf} First of all, by \cref{Prop3}, for an arbitrary $\ep>0$, we can find a sufficiently large $a\geq1-b$ such that $\ep\geq 12M|\bOm|^{-a}$. For $r=\ep/(12M)$, the inequality
\begin{align*}
\f{1}{|\bOm|}\left\|\bu^{\La}-\bsf\right\|_{\ell_2(\bOm)}^2\leq\ep+\f{16}{3}\eta^2
\end{align*}
with probability at least
\begin{align*}
1-\mN\left(\msM,\f{\ep}{12M}\right)\exp\left(-\f{3m\ep}{256M^2}\right)\geq1-\exp\left(\f{480M^2(2a+b) C_{\bsf}\left(d+2C_{\bsf}\right)\left|\bOm\right|^b\log_2\left|\bOm\right|}{\ep}-\f{3m\ep}{256M^2}\right),
\end{align*}
where the last inequality comes from \cref{Th1}. Then, choosing a special $\ep^*$ to be the unique positive solution to
\begin{align}\label{Equation}
\f{480M^2(2a+b)C_{\bsf}\left(d+2C_{\bsf}\right)|\bOm|^b\log_2|\bOm|}{\ep}-\f{3m\ep}{256M^2}=\ln\f{1}{|\bOm|},
\end{align}
we have
\begin{align*}
\f{1}{|\bOm|}\left\|\bu^{\La}-\bsf\right\|_{\ell_2(\bOm)}^2\leq\ep^*+\f{16}{3}\eta^2
\end{align*}
with probability at least $1-|\bOm|^{-1}$. Indeed, solving \cref{Equation} gives
\begin{align}\label{SpecialEpsilon}
\begin{split}
\ep^*&=\f{128M^2}{3m}\left(\ln|\bOm|+\sqrt{\ln^2|\bOm|+\f{45}{2}m(2a+b)C_{\bsf}\left(d+2C_{\bsf}\right)|\bOm|^b\log_2|\bOm|}\right)\\
&\leq\f{64M^2}{3m}\left(4\ln|\bOm|+3\sqrt{10m(2a+b)C_{\bsf}\left(d+2C_{\bsf}\right)|\bOm|^b\log_2|\bOm|}\right)\\
&\leq\f{64}{3}M^2\left(4+3\sqrt{10(2a+b)C_{\bsf}\left(d+2C_{\bsf}\right)}\right)\rho^{-1/2}|\bOm|^{-\f{1-b}{2}}\left(\log_2|\bOm|\right)^{3/2}.
\end{split}
\end{align}
In addition, from \cref{SpecialEpsilon}, we further have
\begin{align*}
\ep^*\geq\f{64M^2}{\sqrt{m}}\sqrt{10(2a+b)C_{\bsf}\left(d+2C_{\bsf}\right)\left|\bOm\right|^b\log_2\left|\bOm\right|}\geq64M\rho^{-1/2}\sqrt{C_{\bsf}\log_2\left|\bOm\right|}\left|\bOm\right|^{-\f{1-b}{2}}\geq12M\left|\bOm\right|^{-a}
\end{align*}
as $a\geq1-b$ with $b\in[0,1)$. This completes the proof.\qquad$\square$
\end{pf}

To conclude this section, we demonstrate the theoretical error bound in \cref{MainResult} and the empirical restoration error under different settings of $|\bOm|$ and $\rho$ through numerical simulations. More precisely, we fix $N=2^J$ and $\eta=0$, and we consider the following noise-free case
\begin{align}\label{NoiseFreeTVInpainting}
\min~\left\|\Na\bu\right\|_1~~~\text{subject to}~~~\bu[\bk]=\bsf[\bk],~~\bk\in\La,
\end{align}
so that \cref{MainResult} is rewritten as
\begin{align}\label{MainResult:NoiseFree}
\f{1}{2^{2J}}\left\|\bu^{\La}-\bsf\right\|_{\ell_2(\bOm)}^2\leq\wt{c}\rho^{-1/2}J^{3/2}2^{-J/2}
\end{align}
with probability at least $1-2^{-2J}$.

\begin{figure}[tp]
\centering
\subfloat[]{\label{PhantomResultsSampleDensity}\includegraphics[width=8.1cm]{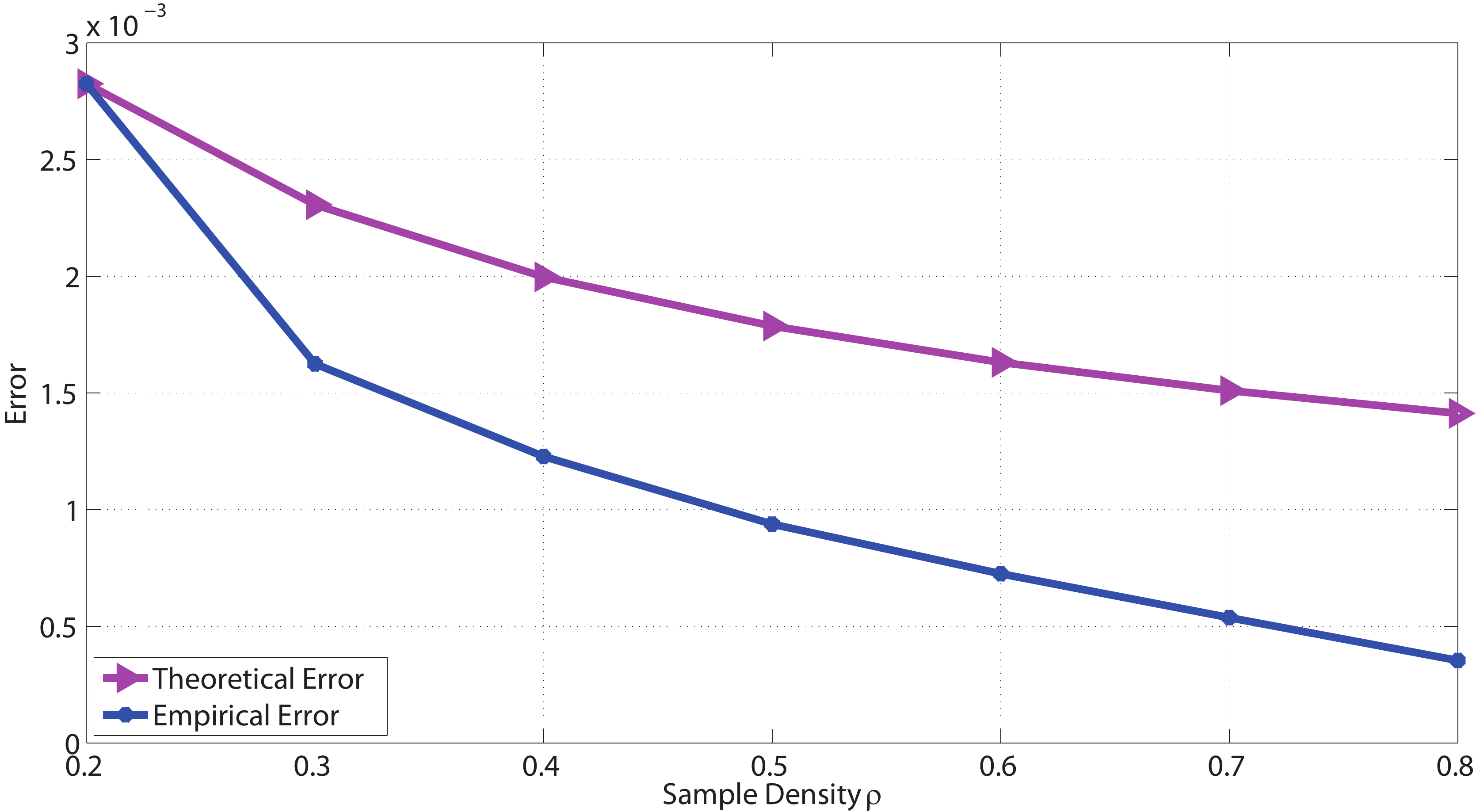}}\hspace{0.001cm}
\subfloat[]{\label{PhantomResultsResolution}\includegraphics[width=8.1cm]{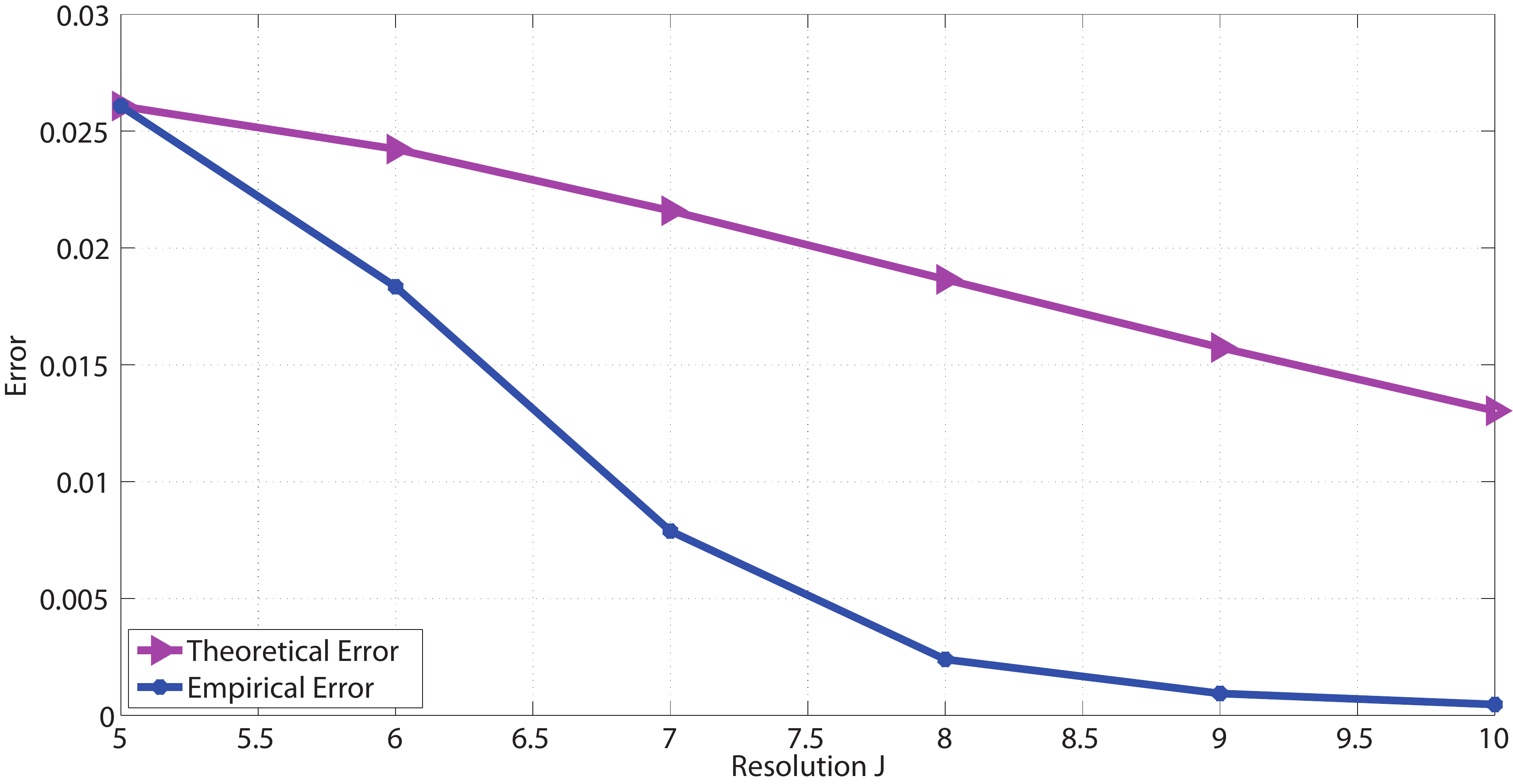}}\hspace{0.001cm}
\caption{\cref{PhantomResultsSampleDensity} describes the simulation results when $N=512$ is fixed, and \cref{PhantomResultsResolution} depicts the simulation results when $\rho=0.5$ is fixed.}\label{PhantomTVInpaintingResultsAnalysis}
\end{figure}

We generate $\bsf$ by using the Matlab built-in function ``phantom($2^J$)''. To see the behavior of error with respect to the sample density $\rho$, we fix $N=512$ (i.e. $J=9$), and for each $\rho\in\left\{0.2,0.3,0.4,0.5,0.6,0.7,0.8\right\}$, we test \cref{NoiseFreeTVInpainting} with $100$ realizations of $\La$. To see the behavior of error with respect to $|\bOm|$, we fix $\rho=0.5$, and for each $J\in\left\{5,6,7,8,9,10\right\}$ (or the resolution $J$), we again test \cref{NoiseFreeTVInpainting} with $100$ realizations of $\La$. In any case, \cref{NoiseFreeTVInpainting} is solved by the split Bregman algorithm (e.g. \cite{T.Goldstein2009}), and we choose the largest empirical error to compare with the theoretical error in \cref{MainResult:NoiseFree}. Since it is in general difficult to determine the constants explicitly, we calculate the above error by assuming that the equality holds in the worst case, i.e. in the lowest sample density, or the lowest resolution case. The results are shown in \cref{PhantomTVInpaintingResultsAnalysis}. Specifically, \cref{PhantomResultsSampleDensity} demonstrates the results when $N=512$ is fixed and $\rho$ is varying, and \cref{PhantomResultsResolution} depicts the results when $\rho=0.5$ is fixed and $J$ is varying. We can easily see that, in any case, the empirical restoration error does not exceed the theoretical error in \cref{MainResult:NoiseFree}, which empirically demonstrates that \cref{Th2} provides a reasonable upper bound for the restoration error with high probability.

\section{Application to sparse gradient restoration}\label{SparseRestoration}

In this section, we connect our main result to the missing data restoration from random samples given that the original data $\bsf$ has a sparse gradient $\Na\bsf$. The error analysis for the sparse gradient restoration have been well established in the literature, and most restoration guarantees are in general based on the context of compressed sensing (e.g. \cite{B.Adcock2021,F.Krahmer2014,Poon2015}). The common concept is that the restoration error is bounded by the so-called $s$-term approximation error of $\Na\bsf$, thereby establishing the restoration guarantee for the $s$-sparse $\Na\bsf$. While these works mainly focus on the Gaussian sensing matrix \cite{J.F.Cai2015}, the Fourier undersampling \cite{B.Adcock2021,F.Krahmer2014,Poon2015}, and the Walsh sampling \cite{B.Adcock2021}, we explore the approximation property of \cref{TVMinimization} with respect to the sparsity of $\Na\bsf$ in terms of \cref{Th2}.

\begin{theorem}\label{Th3} Assume that $\bsg\in[0,M]^{\La}$ in \cref{Linear_IP_Incomplete} is a nonconstant vector which satisfies \cref{Varianceofg}. Let $\left\|\Na\bsf\right\|_0=s$, where $1\leq s\leq|\bOm|-1$, and
\begin{align}\label{Sparsity}
\left\|\Na\bsf\right\|_0=\left|\left\{\bk\in\bOm:\bk+\e_j\in\bOm,~j=1,\cdots,d,~\text{and}~\sum_{j=1}^d\left|\bsf[\bk+\e_j]-\bsf[\bk]\right|\neq0\right\}\right|.
\end{align}
Let $\bu^{\La}$ be a solution to \cref{TVMinimization}. Then the following inequalities
\begin{align}\label{Result-PCSignal2}
\f{1}{|\bOm|}\left\|\bu^{\La}-\bsf\right\|_{\ell_2(\bOm)}^2\leq\wt{c}\rho^{-1/2}\left(\log_2|\bOm|\right)^{3/2}\sqrt{\f{s}{|\bOm|}}+\f{16}{3}\eta^2
\end{align}
and
\begin{align}\label{Result-PCSignal}
\f{1}{|\bOm|}\left\|\bu^{\La}-\bsf\right\|_{\ell_2(\bOm)}^2\leq\wt{c}\left(\log_2|\bOm|\right)^{3/2}\sqrt{\f{s}{m}}+\f{16}{3}\eta^2
\end{align}
hold with probability at least $1-|\bOm|^{-1}$. In \cref{Result-PCSignal2,Result-PCSignal}, $\wt{c}$ is defined as
\begin{align*}
\wt{c}=\f{128}{3}M^2\left(2+3\sqrt{5\left(2a+1\right)M\left(d+4M\right)}\right)
\end{align*}
for some $a\geq1$.
\end{theorem}

\begin{pf} Notice that we have
\begin{align*}
\left\|\Na\bsf\right\|_1\leq s\left\|\Na\bsf\right\|_{\infty}=\left\|\Na\bsf\right\|_{\infty}|\bOm|^{\log_{|\bOm|}s},
\end{align*}
where
\begin{align*}
\left\|\Na\bsf\right\|_{\infty}=\max_{j=1,\cdots,d}\left\{\left|\bsf[\bk+\e_j]-\bsf[\bk]\right|:\bk,\bk+\e_j\in\bOm\right\}.
\end{align*}
Then for each $j=1,\cdots,d$, we have
\begin{align*}
\left|\bsf[\bk+\e_j]-\bsf[\bk]\right|\leq2\left\|\bsf\right\|_{\ell_{\infty}(\bOm)}\leq2M.
\end{align*}
Choose $C_{\bsf}=2M$ and $b=\log_{|\bOm|}s$. For $a\geq1$, we have $a\geq1-\log_{|\bOm|}s$, so \cref{MainResult} becomes
\begin{align*}
\f{1}{|\bOm|}\left\|\bu^{\La}-\bsf\right\|_{\ell_2(\bOm)}^2&\leq\f{128}{3}M^2\left(2+3\sqrt{5\left(2a+\log_{|\bOm|}s\right)M\left(d+4M\right)}\right)\rho^{-1/2}\left(\log_2|\bOm|\right)^{3/2}\sqrt{\f{s}{|\bOm|}}+\f{16}{3}\eta^2\\
&\leq\f{128}{3}M^2\left(2+3\sqrt{5\left(2a+1\right)M\left(d+4M\right)}\right)\rho^{-1/2}\left(\log_2|\bOm|\right)^{3/2}\sqrt{\f{s}{|\bOm|}}+\f{16}{3}\eta^2.
\end{align*}
Using the fact that $\rho=m/|\bOm|$, we further have
\begin{align*}
\f{1}{|\bOm|}\left\|\bu^{\La}-\bsf\right\|_{\ell_2(\bOm)}^2\leq\f{128}{3}M^2\left(2+3\sqrt{5\left(2a+1\right)M\left(d+4M\right)}\right)\left(\log_2|\bOm|\right)^{3/2}\sqrt{\f{s}{m}}+\f{16}{3}\eta^2,
\end{align*}
and this completes the proof.\qquad$\square$
\end{pf}

\cref{Th3} tells us that for a fixed $\bOm$, when we solve the TV minimization for the noise-free setting:
\begin{align}\label{Minimizer_EqualConstraint}
\bu^{\La}=\argmin_{\bu}\left\{\left\|\Na\bu\right\|_1:\bu[\bk]=\bsf[\bk],~\bk\in\La\right\},
\end{align}
the error satisfies
\begin{align*}
\f{1}{|\bOm|}\left\|\bu^{\La}-\bsf\right\|_{\ell_2(\bOm)}^2=O\left(\sqrt{\f{s}{m}}\right)
\end{align*}
with high probability. We would like to mention that this error bound cannot be equal to $0$ to guarantee the exact restoration even in the noise-free setting. Notice that, from the viewpoint of compressed sensing, the sensing matrix of our setting will be $\bmR_{\La}$, where $\bmR_{\La}\bsf[\bk]=\bsf[\bk]$ for $\bk\in\La$, and $\bmR_{\La}\bsf[\bk]=0$ for $\bk\in\bOm\setminus\La$. Since $\La$ is randomly chosen from the uniform distribution of $\bOm$, the sensing matrix $\bmR_{\La}$ may not satisfy the concentration inequality in \cite{H.Rauhut2008}, which means that the measurement $\bsg$ may not contain sufficient information for the exact restoration \cite{J.F.Cai2011}.

\begin{example} For a better explanation, we consider the one dimensional case ($d=1$). Let $N=2N_0$ for some $N_0\in\N$, and let $\bsf\in\R^N$ be defined as
\begin{align*}
\bsf[k]=\left\{\begin{array}{cll}
0~&\text{if}&k=0,\cdots,N_0-1\vspace{0.25em}\\
1~&\text{if}&k=N_0,\cdots,N-1.
\end{array}\right.
\end{align*}
Obviously, we have $s=1$. Assume that $\La$ is uniformly and randomly drawn from all $m$-subsets of $\bOm$. We claim that there exists a solution $\bu^{\La}$ in \cref{Minimizer_EqualConstraint} such that
\begin{align*}
\f{1}{|\bOm|}\left\|\bu^{\La}-\bsf\right\|_{\ell_2(\bOm)}^2\geq\f{1}{N}
\end{align*}
with probability at least $1-m/N$. To see this, define $k_1$ and $k_2$ as
\begin{align*}
k_1&=\max\left\{k\in\bOm:k\in\La\cap\left\{0,\cdots,N_0-1\right\}\right\}\\
k_2&=\min\left\{k\in\bOm:k\in\La\cap\left\{N_0+1,\cdots,N-1\right\}\right\}.
\end{align*}
For each $L\in\left\{k_1,k_1+1,\cdots,k_2-1,k_2\right\}$, we introduce
\begin{align}\label{AllPossibleMinimizers}
\bu^{\La,L}=\left\{\begin{array}{cll}
0~&\text{if}&k=0,\cdots,L-1\vspace{0.25em}\\
1~&\text{if}&k=L,\cdots,N-1.
\end{array}\right.
\end{align}
Then obviously, $\bu^{\La,L}$ is a solution of \cref{Minimizer_EqualConstraint} for each $L=k_1,\cdots,k_2$, and
\begin{align*}
\f{1}{|\bOm|}\left\|\bu^{\La,k_1}-\bsf\right\|_{\ell_2(\bOm)}^2&=\f{N_0-1-k_1}{N}\\
\f{1}{|\bOm|}\left\|\bu^{\La,k_2}-\bsf\right\|_{\ell_2(\bOm)}^2&=\f{k_2-N_0}{N}.
\end{align*}
Hence, if $\La\cap\left\{N_0\right\}=\emptyset$, there exists a solution $\bu^{\La}$ to \cref{Minimizer_EqualConstraint} such that
\begin{align*}
\f{1}{|\bOm|}\left\|\bu^{\La}-\bsf\right\|_{\ell_2(\bOm)}^2\geq\f{1}{N}.
\end{align*}
In addition,
\begin{align*}
\PP\left\{\La\cap\left\{N_0\right\}=\emptyset\right\}=\binom{N-1}{m}/\binom{N}{m}=\f{N-m}{N}=1-\f{m}{N},
\end{align*}
and this completes the proof.
\end{example}

\section{Application to two dimensional $BV$ function approximation}\label{TVInpainting}

This section is devoted to the connection of the total variation image inpainting (e.g. \cite{T.F.Chan2001/02}) to the underlying function approximation. In the literature, there are various numerical algorithms for the total variation minimization in \cite{A.Chambolle2011,T.Goldstein2009,Y.Wang2008,C.Wu2010,X.Zhang2011} with a guaranteed convergence to the minimizer. Hence, we are able to analyze the approximation property of these numerical algorithms. In addition, based on the finite element approximation of $BV$ functions in \cite{S.Bartels2015}, we connect the error analysis in the discrete setting to the approximation of underlying $BV$ function from which a discrete image comes. In what follows, we restrict our discussions for the real-valued function of two variables ($d=2$), as the images can be treated as discrete samples of two variable functions \cite{J.F.Cai2011}. Note, however, that for more general multivariate functions, the discussions are almost the same with a slight modification.

All functions we consider are defined on the square domain $\Om=[0,1)^2$, and we assume for simplicity that $\bOm$ is a $2^J\times 2^J$ cartesian grid defined as
\begin{align*}
\bOm=\left\{0,1,\cdots,2^J-1\right\}^2\simeq2^{-J}\Z^2\cap\Om.
\end{align*}
In other words, we implicitly identify a $2^J\times2^J$ grid $\bOm$ with a $2^J\times2^J$ discrete mesh of $\Om$. Note, however, that it is not difficult to extend our discussion to the generic regular square grid $\bOm\simeq N^{-1}\Z^2\cap\Om$. To establish a suitable approximation analysis, we assume that the functions on $\Om$ with the ones on $\R^2$ with fundamental period of each variable to be $1$.

Recall that a function $u\in BV(\Om)$ (is of bounded variation) if $u\in L_1(\Om)$ and its distributional first order derivative is a Radon measure. To simplify the notation, we use $\na u$ to denote such a measure. We define the total variation (TV) of $u\in BV(\Om)$ by
\begin{align}\label{FunctionTV}
\TV(u)=\int_{\Om}\rd\left\|\na u\right\|_1
\end{align}
with $\left\|\cdot\right\|_1$ being the $\ell_1$ norm in $\R^2$. Notice that the above TV is the anisotropic TV, which is the usual choice in the study of multidimensional nonlinear conservation laws (e.g. \cite{S.Bartels2014}).

To begin with, let $\varphi=1_{\Om}$. We assume that, for each $\bk\in\bOm$, the discrete samples $\bsf$ are obtained via
\begin{align}\label{Sampling1}
\bsf[\bk]=2^{J}\left\la f,\varphi_{J,\bk}\right\ra,
\end{align}
where $\varphi_{J,\bk}=2^{J}\varphi(2^J\cdot-\bk)$. In other words, we assume that the discrete samples are obtained by taking the local averages of the underlying function on the square $\bfQ_{\bk}:=\prod_{j=1}^2[2^{-J}k_j,2^{-J}(k_j+1))$.

Given the discrete samples $\bsf$ on $\bOm$, we use the following interpolated function
\begin{align}\label{Interpolation1}
f_J=\sum_{\bk\in\bOm}\bsf[\bk]\phi(2^J\cdot-\bk)=\sum_{\bk\in\bOm}\left\la f,\varphi_{J,\bk}\right\ra\phi_{J,\bk}
\end{align}
to approximate $f$. In \cref{Interpolation1}, $\phi(\x)=B_2(x_1)B_2(x_2)$ is a tensor product piecewise linear B-spline:
\begin{align*}
B_2(x)=\max\left\{1-|x|,0\right\},
\end{align*}
and we implicitly identify $\phi(2^J\cdot-\bk)$ with its periodized version
\begin{align*}
\phi^{\mathrm{per}}(2^J\cdot-\bk)=\sum_{\bk'\in\Z^2}\phi(2^J\cdot-\bk+2^J\bk')
\end{align*}
with a slight abuse of notation.

In the literature, there are extensive studies on the approximation order of the interpolated function to the underlying function. Most of them are related with the property of the basis function $\phi$, and require a high order regularity of an underlying function $f$. Briefly speaking, if $\phi$ satisfies the Strang-Fix condition of a certain order and its Fourier transform $\wh{\phi}$ is such that $1-\left|\wh{\phi}(\oom)\right|^2$ has the same order of zeroes at $\oom=\0$, then for a sufficiently smooth $f$, the interpolated function $f_J$ has the approximation of this order to $f$ \cite{I.Daubechies2003,M.J.Johnson2009}. Indeed, for the following harmonic inpainting (e.g. \cite{T.F.Chan2001/02})
\begin{align*}
\min~\left\|\Na\bu\right\|_2^2~~~\text{subject to}~~~\f{1}{|\La|}\sum_{\bk\in\La}\left|\bu[\bk]-\bsg[\bk]\right|^2\leq\eta^2,
\end{align*}
the asymptotic approximation analysis can be done similarly to \cite{M.J.Johnson2009}. Unlike the aforementioned ideas which requires a high regularity of $f$, we only assume that $f$ is of bounded variation (i.e. $f\in BV(\Om)$) to ensure $\left\|\Na\bsf\right\|_1<\infty$. In fact, using the finite element approximation \cite{S.Bartels2015}, we present the approximation order of $f_J$ to $f$, as given in \cref{Th5}. The proof is postponed to \ref{PfTh5}.

\begin{theorem}\label{Th5} Assume that $f\in BV(\Om)\cap L_{\infty}(\Om)$. For $J\in\N$, let $f_J$ be defined as \cref{Interpolation1} with $\bsf$ in \cref{Sampling1}. Then we have
\begin{align}\label{InterpolationApproximation}
\left\|f_J-f\right\|_{L_2(\Om)}^2\leq\left(16+4\pi^{1/2}\right)2^{-J}\TV(f)\left\|f\right\|_{L_{\infty}(\Om)}.
\end{align}
\end{theorem}

Combining \cref{Th2,Th5}, we are able to present \cref{Th4} to connect the solution to the discrete problem \cref{TVMinimization} to the underlying function approximation. Briefly speaking, as long as the mesh is sufficiently dense, we have a good opportunity to obtain a reasonable approximation of the underlying true image $\bsf$ by solving \cref{TVMinimization}. Moreover, the interpolation of the restored image gives a good approximation of the original function where the discrete image comes from, with the high probability. The proof is in \cref{PfTh4}.

\begin{theorem}\label{Th4} Assume that $f\in BV(\Om)\cap L_{\infty}(\Om,[0,M])$ is not identically constant. Let $\bu^{\La}$ be a solution to \cref{TVMinimization} with $\bsg$ in \cref{Linear_IP_Incomplete} generated by $\bsf$ in \cref{Sampling1}. Then the inequality
\begin{align}\label{MainResult:Recall}
\f{1}{2^{2J}}\left\|\bu^{\La}-\bsf\right\|_{\ell_2(\bOm)}^2\leq\wt{c}\rho^{-1/2}J^{3/2}2^{-J/2}+\f{16}{3}\eta^2
\end{align}
holds with probability at least $1-2^{-2J}$, with a constant $\wt{c}$ independent of $J$, $\rho$, and $\eta$. Moreover, let $u_J^{\La}$ be defined as
\begin{align}\label{Interpolation2}
u_J^{\La}=\sum_{\bk\in\bOm}\bu^{\La}[\bk]\phi(2^J\cdot-\bk).
\end{align}
Then the following inequality
\begin{align}\label{Coro1:Result}
\left\|u_J^{\La}-f\right\|_{L_2(\Om)}^2\leq C_1\rho^{-1/2}J^{3/2}2^{-J/2}+C_2\eta^2+C_32^{-J}
\end{align}
also holds with probability at least $1-2^{-2J}$, where $C_1$, $C_2$, and $C_3$ are independent of $J$, $\rho$, and $\eta$.
\end{theorem}

\begin{rmk} Note that, if $J\in\N$ satisfies $J\geq-2\log_2\eta$, \cref{Coro1:Result} becomes
\begin{align*}
\left\|u_J^{\La}-f\right\|_{L_2(\Om)}^2\leq C_1\rho^{-1/2}J^{3/2}2^{-J/2}+\left(C_2+C_3\right)\eta^2.
\end{align*}
This further means that, if the mesh is sufficiently dense (i.e. $J$ is sufficiently large), the $L_2$ distance between the interpolated function $u_J^{\La}$ and the original underlying function $f$ becomes bounded by the restoration error of the discrete image restoration problem \cref{TVMinimization} only.
\end{rmk}

To conclude this section, we further discuss the approximation of piecewise constant function from the discrete sparse gradient restoration problem. To be more precise, let $f$ be defined as
\begin{align}\label{fModel}
f(\x)=\sum_{l=1}^L\alpha_l1_{\Om_l}(\x)
\end{align}
where $\alpha_l\in\R$, $\Om_l\subseteq\Om$, and \cref{fModel} is expressed with the smallest number of characteristic functions such that $\Om_j$'s are pairwise disjoint. Let $\bsf$ be defined as \cref{Sampling1}. From \cref{Sampling1}, we have
\begin{align*}
\bsf[\bk+\e_j]-\bsf[\bk]=2^J\left\la f,\varphi_{J,\bk+\e_j}\right\ra-2^J\left\la f,\varphi_{J,\bk}\right\ra=2^J\left(\left\la f,\varphi_{J,\bk+\e_j}-\varphi_{J,\bk}\right\ra\right).
\end{align*}
Noting that
\begin{align*}
1_{[0,1)}(\cdot-1)-1_{[0,1)}=-\f{\rd}{\rd x}B_2(\cdot-1),
\end{align*}
the direct computations gives
\begin{align*}
\bsf[\bk+\e_j]-\bsf[\bk]=\left\la\f{\p f}{\p x_j},\wt{\varphi}_{j,J,\bk+\e_j}\right\ra,
\end{align*}
where $\wt{\varphi}_j$ is defined as
\begin{align*}
\wt{\varphi}_1(\x)=B_2(x_1)1_{[0,1)}(x_2)~~~\text{and}~~~\wt{\varphi}_2(\x)=1_{[0,1)}(x_1)B_2(x_2).
\end{align*}
Therefore, $\bsf[\bk+\e_j]-\bsf[\bk]=0$ if and only if
\begin{align*}
\su\left(\wt{\varphi}_{j,J,\bk+\e_j}\right)=\overline{\bfQ}_{\bk}\cup\overline{\bfQ}_{\bk+\e_j}\subseteq\Om_l^{\circ}~~~\text{for some}~~l=1,\cdots,L.
\end{align*}
For each $l=1,\cdots,L$, we define
\begin{align*}
\I_l=\left\{\bk\in\bOm:\overline{\bfQ}_{\bk}\cup\overline{\bfQ}_{\bk+\e_1}\cup\overline{\bfQ}_{\bk+\e_2}\subseteq\Om_l^{\circ}\right\},
\end{align*}
and we let
\begin{align*}
\I=\bigcup_{l=1}^L\I_l:=\left\{\bk\in\bOm:\overline{\bfQ}_{\bk}\cup\overline{\bfQ}_{\bk+\e_1}\cup\overline{\bfQ}_{\bk+\e_2}\subseteq\Om_l^{\circ}~~\text{for some}~~l=1,\cdots,L-1\right\}.
\end{align*}
Denote $\bS:=\bOm\setminus\I$. Obviously, we have $\left\|\Na\bsf\right\|_0=\left|\bS\right|$, and by applying \cref{Th3,Th5}, we obtain
\begin{align}\label{Result-PCSignal:Recall}
\f{1}{2^{2J}}\left\|\bu^{\La}-\bsf\right\|_{\ell_2(\bOm)}^2\leq\wt{c}\rho^{-1/2}J^{3/2}2^{-J}\left|\bS\right|^{1/2}+\f{16}{3}\eta^2
\end{align}
and
\begin{align}\label{Coro1:Result:Recall}
\left\|u_J^{\La}-f\right\|_{L_2(\Om)}^2\leq C_1\rho^{-1/2}J^{3/2}2^{-J}\left|\bS\right|^{1/2}+C_2\eta^2+C_32^{-J},
\end{align}
with probability at least $1-2^{-2J}$ where constants $\wt{c}$, $C_1$, $C_2$, and $C_3$ are all independent of $J$, $\rho$, and $\eta$. However, it should be noted that $\left|\bS\right|=\left\|\Na\bsf\right\|_0$ may not necessarily satisfy $\left|\bS\right|\leq2^{2bJ}$ with $b<1/2$ as it is related to the geometry of edges $\bigcup_{l=1}^L\p\Om_l$. Hence, the above estimates \cref{Result-PCSignal:Recall,Coro1:Result:Recall} will be worse than \cref{Th4}.

\section{Technical proofs}\label{TechnicalProofs}

This section is devoted to the technical details left in the previous sections. Mainly, we focus on the proof of \cref{Th1,Th5,Th4}.

\subsection{Proof of \cref{Th1}}\label{PfTh1}

\cref{Th1} is to estimate the covering number. The proof follows the line similar to \cite[Theorem 2.4]{J.F.Cai2011}. However, since we improve \cite[Theorem 2.4]{J.F.Cai2011} by relaxing the constraint of the radius $r$, we include the detailed proof for the sake of completeness. Notice that it is obvious that $\msM\subseteq\wt{\msM}$ and $\msN\left(\msM,r\right)\leq\msN\left(\wt{\msM},r\right)$ where $\wt{\msM}$ is defined as \cref{HypothesisSpaceRelax}. Hence, it suffices to bound the covering number $\msN\left(\wt{\msM},r\right)$. In addition, it is easy to see that if there exists a finite set $F\subseteq\wt{\msM}$ such that
\begin{align*}
\wt{\msM}\subseteq\bigcup_{\bq\in F}\left\{\bu:\left\|\bu-\bq\right\|_{\ell_{\infty}(\bOm)}\leq r\right\},
\end{align*}
we have $\msN\left(\wt{\msM},r\right)\leq\left|F\right|$. What we need now is to construct an appropriate set $F$ by exploiting the specific structure of $\wt{\msM}$, so that $\left|F\right|$ has an appropriate upper bound.

For this purpose, let $\kappa=\lceil2M/r\rceil$, and we define
\begin{align}\label{QuantizedRange}
\msR=\left\{-\kappa r/2,(-\kappa+1)r/2,\cdots,\kappa r/2\right\}.
\end{align}
By \cite[Lemma 4.4]{J.F.Cai2011}, for each $\bu\in\wt{\msM}$, there exists $Q(\bu)\in\msR^{\bOm}$ such that
\begin{align*}
\left\|\bu-Q(\bu)\right\|_{\ell_{\infty}(\bOm)}\leq r/2~~~\text{and}~~~\left\|\Na\left(Q(\bu)\right)\right\|_1\leq\left\|\Na\bu\right\|_1.
\end{align*}
Let
\begin{align*}
\wt{F}=\left\{\bq\in\ell_{\infty}(\bOm):\bq=Q(\bu)~~\text{for some}~~\bu\in\wt{\msM}\right\}\subseteq\msR^{\bOm}.
\end{align*}
Notice that for each $\bq\in\wt{F}$, there may be more than one $\bu\in\wt{\msM}$ such that $\bq=Q(\bu)$.

For each $\bq\in\wt{F}$, choose $\bu_{\bq}\in\wt{\msM}$ such that $\left\|\bu_{\bq}-\bq\right\|_{\ell_{\infty}(\bOm)}\leq r/2$, and define $F=\left\{\bu_{\bq}:\bq\in\wt{F}\right\}$. For an arbitrary $\bu\in\wt{\msM}$, there exists $\bq\in\wt{F}$ such that $\left\|\bu-\bq\right\|_{\ell_{\infty}(\bOm)}\leq r/2$. This implies
\begin{align*}
\left\|\bu-\bu_{\bq}\right\|_{\ell_{\infty}(\bOm)}\leq\left\|\bu-\bq\right\|_{\ell_{\infty}(\bOm)}+\left\|\bq-\bu_{\bq}\right\|_{\ell_{\infty}(\bOm)}\leq r,
\end{align*}
by the definition of $\bu_{\bq}$. Therefore,
\begin{align*}
\wt{\msM}\subseteq\bigcup_{\bu_{\bq}\in F}\left\{\bu:\left\|\bu-\bu_{\bq}\right\|_{\ell_{\infty}(\bOm)}\leq r\right\}~~\text{and}~~\msN\left(\wt{\msM},r\right)\leq\left|F\right|\leq\left|\wt{F}\right|.
\end{align*}
Thus, the covering number $\msN\left(\wt{\msM},r\right)$ is bounded by any upper bound of $\left|\wt{F}\right|$. Notice that each $\bq\in\wt{F}$ is uniquely determined by $\Na\bq$ and $\bq[1,\cdots,1]$. Since $\wt{F}$ is a subset of $\msR^{\bOm}$, there are $2\kappa+1$ choices for $\bq[1,\cdots,1]$. It remains to count the number of choices in $\Na\bq$. Define
\begin{align*}
\Na\wt{F}=\left\{\Na\bq:\bq\in\wt{F}\right\}.
\end{align*}
Then we need to bound $\left|\Na\wt{F}\right|$.

To do this, we first consider the uniform upper bound of $\left\|\Na\bq\right\|_1$ for $\bq\in\wt{F}$. By the definition of $\wt{F}$ and \cite[Lemma 4.4]{J.F.Cai2011}, for each $\bq\in\wt{F}$, there exists $\bu\in\wt{\msM}$ such that $\left\|\Na\bq\right\|_1\leq\left\|\Na\bu\right\|_1$. Since $\left\|\Na\bsf\right\|_1\leq C_{\bsf}\left|\bOm\right|^b$ for some $b\in[0,1]$ by assumption, we further have
\begin{align*}
\left\|\Na\bq\right\|_1\leq\left\|\Na\bu\right\|_1\leq\left\|\Na\bsf\right\|_1\leq C_{\bsf}\left|\bOm\right|^b.
\end{align*}
Hence, for all $\bq\in\wt{F}$, we have $\left\|\Na\bq\right\|_1\leq Kr/2$, where
\begin{align*}
K=\left\lceil\f{2C_{\bsf}\left|\bOm\right|^b}{r}\right\rceil.
\end{align*}
In addition, since $\bq\in\msR^{\bOm}$, each element of $\Na\bq$ has to be a multiple of $r/2$, which means that the range of $\Na\bq$ is a subset of $\left\{-Kr/2,-(K-1)r/2,\cdots,Kr/2\right\}^d$. Recall that there are $R=d\left(\left|\bOm\right|-\left|\bOm\right|^{(d-1)/d}\right)$ elements in $\Na\bq$. Hence, the bound of $\left|\Na\wt{F}\right|$ can be estimated by the number of possible integer solutions of the following inequality
\begin{align*}
\left|x_1\right|+\left|x_2\right|+\cdots+\left|x_{R-1}\right|+\left|x_R\right|\leq K.
\end{align*}
That is,
\begin{align*}
\left|\Na\wt{F}\right|&\leq1+\sum_{k=1}^K\sum_{s=1}^{\min\{k,R\}}2^s\binom{R}{R-s}\binom{k-1}{s-1}\\
&\leq1+\sum_{k=1}^K2^k\sum_{s=1}^{\min\{k,R\}}\binom{R}{R-s}\binom{k-1}{s-1}=1+\sum_{k=1}^K2^k\binom{R+k-1}{R-1}\\
&\leq\binom{R+K-1}{R-1}\sum_{k=0}^K2^k\leq2^{K+1}\binom{R+K-1}{K}\leq2\left[2\left(R+K-1\right)\right]^K.
\end{align*}
Hence, we have
\begin{align*}
\left|\wt{F}\right|\leq\left(4\kappa+2\right)\left[2\left(R+K-1\right)\right]^K.
\end{align*}
In other words, using $r\geq\left|\bOm\right|^{-a}$ with $a\geq1-b$ and $K-1\leq2C_{\bsf}|\bOm|^{a+b}$, we have
\begin{align*}
\ln\mN\left(\msM,r\right)&\leq\f{2C_{\bsf}\left|\bOm\right|^b}{r}\ln\left(2\left(R+K-1\right)\right)+\ln\left(\f{8M}{r}+2\right)\\
&\leq\f{2C_{\bsf}\left|\bOm\right|^b}{r}\left[\ln\left(2d\left|\bOm\right|+4C_{\bsf}\left|\bOm\right|^{a+b}\right)+\ln\f{10M}{r}\right]\leq\f{2C_{\bsf}\left|\bOm\right|^b}{r}\ln\left(\left(d+2C_{\bsf}\right)20M\left|\bOm\right|^{2a+b}\right)
\end{align*}
where we use the fact that $a+b\geq1$ from the choice of $a$ and $b$ in the final inequality. Therefore, we have
\begin{align*}
\ln\mN\left(\msM,r\right)\leq\f{40(2a+b)MC_{\bsf}\left(d+2C_{\bsf}\right)\left|\bOm\right|^b}{r}\log_2\left|\bOm\right|,
\end{align*}
and this completes the proof.

\subsection{Proof of \cref{Th5}}\label{PfTh5}

In this section, we prove \cref{Th5}. To begin with, notice that, by the standard density argument of $BV(\Om)$ (e.g. \cite{H.Attouch2014}): for each $f\in BV(\Om)$, there exists $f_n\in C^{\infty}(\Om)$ such that
\begin{align*}
\lim_{n\to\infty}f_n=f~~~\text{in}~~L_1(\Om)~~\text{and}~~\lim_{n\to\infty}\int_{\Om}\left\|\na f_n(\x)\right\|_1\rd\x=\TV(f),
\end{align*}
it suffices to prove \cref{InterpolationApproximation} for $f\in C^{\infty}(\Om)$ (i.e. $\na f$ is defined in the classical sense). Then since the constant in \cref{InterpolationApproximation} is independent of the choice of $f$, \cref{InterpolationApproximation} holds for $f\in BV(\Om)$ as well.

First of all, by the interpolation of $L_p$ spaces (e.g. \cite{Folland1999}), it suffices to estimate $\left\|f_J-f\right\|_{L_p(\Om)}$ for $p=1$ and $p=\infty$, as we then have
\begin{align*}
\left\|f_J-f\right\|_{L_2(\Om)}\leq\left\|f_J-f\right\|_{L_1(\Om)}^{1/2}\left\|f_J-f\right\|_{L_{\infty}(\Om)}^{1/2}.
\end{align*}
For $p=\infty$, since $0\leq\phi(2^J\cdot-\bk)\leq1$ and it forms a partition of unity, we have
\begin{align}\label{LinfinityError}
\left\|f_J-f\right\|_{L_{\infty}(\Om)}\leq\left\|f_J\right\|_{L_{\infty}(\Om)}+\left\|f\right\|_{L_{\infty}(\Om)}\leq\left\|\bsf\right\|_{\ell_{\infty}(\bOm)}+\left\|f\right\|_{L_{\infty}(\Om)}\leq2\left\|f\right\|_{L_{\infty}(\Om)}
\end{align}
where we used the H\"older's inequality \cite{Folland1999} in the last inequality.

For $p=1$, we note that
\begin{align*}
\left\|f_J-f\right\|_{L_1(\Om)}&=\sum_{\bk\in\bOm}\int_{\bfQ_{\bk}}\left|f_J(\x)-f(\x)\right|\rd\x\\
&\leq\sum_{\bk\in\bOm}\int_{\bfQ_{\bk}}\left|f_J(\x)-\bsf[\bk]\right|\rd\x+\sum_{\bk\in\bOm}\int_{\bfQ_{\bk}}\left|\bsf[\bk]-f(\x)\right|\rd\x.
\end{align*}
For the first term, let $\x\in\bfQ_{\bk}$. Since $0\leq\phi(2^J\cdot-\bk)\leq1$ and it forms a partition of unity, we have
\begin{align*}
\left|f_J(\x)-\bsf[\bk]\right|=\left|\sum_{\bsl\in\msT_{\bk}}\bsf[\bsl]\phi(2^J\x-\bsl)-\bsf[\bk]\sum_{\bsl\in\msT_{\bk}}\phi(2^J\x-\bsl)\right|\leq\sum_{\bsl\in\msT_{\bk}}\left|\bsf[\bsl]-\bsf[\bk]\right|\phi(2^J\x-\bsl),
\end{align*}
where $\msT_{\bk}=\left\{\bk,\bk+\e_1,\bk+\e_2,\bk+\e_1+\e_2\right\}$ denotes the vertices of $\bfQ_{\bk}$. Notice that we have
\begin{align*}
1_{[0,1)}(\cdot-1)-1_{[0,1)}=-\f{\rd}{\rd x}B_2(\cdot-1).
\end{align*}
Hence, from \cref{Sampling1}, the direct computation gives
\begin{align}\label{DiffvsDerivative}
\bsf[\bk+\e_j]-\bsf[\bk]=\left\la\f{\p f}{\p x_j},\wt{\varphi}_{j,J,\bk+\e_j}\right\ra,
\end{align}
and
\begin{align}\label{DiffvsDerivative2}
\bsf[\bk+\e_1+\e_2]-\bsf[\bk]=\left\la\f{\p f}{\p x_2},\wt{\varphi}_{2,J,\bk+\e_1+\e_2}\right\ra+\left\la\f{\p f}{\p x_1},\wt{\varphi}_{1,J,\bk+\e_1}\right\ra,
\end{align}
where $\wt{\varphi}_j$ is defined as
\begin{align*}
\wt{\varphi}_1(\x)=B_2(x_1)1_{[0,1)}(x_2)~~~\text{and}~~~\wt{\varphi}_2(\x)=1_{[0,1)}(x_1)B_2(x_2).
\end{align*}
From \cref{DiffvsDerivative,DiffvsDerivative2}, we have
\begin{align*}
\sum_{\bsl\in\msT_{\bk}}\left|\bsf[\bsl]-\bsf[\bk]\right|\leq2^{J+1}\int_{\bS_{\bk+\e_1+\e_2}}\rd\left\|\na f\right\|_1
\end{align*}
where $\bS_{\bk}=\prod_{j=1}^2[2^{-J}(k_j-1),2^{-J}(k_j+1)$. Together with the fact that
\begin{align*}
\int_{\bfQ_{\bk}}\sum_{\bsl\in\msT_{\bk}}\phi(2^J\x-\bsl)\rd\x=\int_{\R^2}\phi(2^J\x-\bsl)\rd\x=2^{-2J}~~\text{and}~~\sum_{\bk\in\bOm}1_{\bS_{\bk}}=4,
\end{align*}
we have
\begin{align}\label{L1Error1}
\sum_{\bk\in\bOm}\int_{\bfQ_{\bk}}\left|f_J(\x)-\bsf[\bk]\right|\rd\x\leq2^{1-J}\sum_{\bk\in\bOm}\int_{\bS_{\bk+\e_1+\e_2}}\rd\left\|\na f\right\|_1=2^{3-J}\TV(f).
\end{align}
For the second term, we note that for each $\bk\in\bOm$,
\begin{align*}
\bsf[\bk]=2^J\left\la f,\varphi_{J,\bk}\right\ra=2^{2J}\int_{\bfQ_{\bk}}f(\y)\rd\y.
\end{align*}
Hence, the estimation follows the similar line to \cite[Lemma 7.16]{D.Gilbarg2001}. More precisely, for $\x,\y\in\bfQ_{\bk}$, we have
\begin{align*}
f(\x)-f(\y)=-\int_0^{|\x-\y|}\p_rf(\x+r\ssi)\rd r,~~~\text{where}~~~\ssi=\f{\y-\x}{|\y-\x|}
\end{align*}
and $\p_rf(\x+r\ssi)=\na f(\x+r\ssi)\cdot\ssi$. Integrating over $\bfQ_{\bk}$ with respect to $\y$, we have
\begin{align*}
2^{-2J}\left(f(\x)-\bsf[\bk]\right)=-\int_{\bfQ_{\bk}}\int_0^{|\x-\y|}\p_rf(\x+r\ssi)\rd r\rd\y.
\end{align*}
For the notational simplicity, we introduce $g_{\bk}(\x)$ as
\begin{align}\label{g_k}
g_{\bk}(\x)=\left\{\begin{array}{ccl}
|\p_rf(\x)|~&\text{if}&~\x\in\bfQ_{\bk}\vspace{0.25em}\\
0~&\text{if}&~\x\notin\bfQ_{\bk}.
\end{array}\right.
\end{align}
Then using the polar coordinate, we have
\begin{align*}
\left|f(\x)-\bsf[\bk]\right|&\leq2^{2J}\int_{|\x-\y|<2^{-J+1/2}}\int_0^{\infty}g_{\bk}(\x+r\ssi)\rd r\rd\y\\
&=2^{2J}\int_0^{\infty}\int_{|\ssi|=1}\int_0^{2^{-J+1/2}}g_{\bk}(\x+r\ssi)\rho\rd\rho\rd\ssi\rd r=\int_0^{\infty}\int_{|\ssi|=1}g_{\bk}(\x+r\ssi)\rd\ssi\rd r.
\end{align*}
Let $\z=\x+r\ssi$. Then $r=|\x-\z|$, and from the definition \cref{g_k} of $g_{\bk}$, we have
\begin{align*}
\left|f(\x)-\bsf[\bk]\right|\leq\int_{\bfQ_{\bk}}\f{\left|\na f(\z)\right|}{|\x-\z|}\rd\z\leq\int_{\bfQ_{\bk}}\f{\left\|\na f(\z)\right\|_1}{|\x-\z|}\rd\z.
\end{align*}
where we emphasize the variable of integration for the sake of clarity. By \cite[Lemma 7.12]{D.Gilbarg2001}, we have
\begin{align}\label{L1Error2}
\sum_{\bk\in\bOm}\int_{\bfQ_{\bk}}\left|f(\x)-\bsf[\bk]\right|\rd\x\leq\pi^{1/2}2^{-J+1}\sum_{\bk\in\bOm}\int_{\bfQ_{\bk}}\rd\left\|\na f\right\|_1=\pi^{1/2}2^{-J+1}\TV(f).
\end{align}
Hence, by \cref{L1Error1,L1Error2}, we have
\begin{align}\label{L1Error}
\left\|f_J-f\right\|_{L_1(\Om)}\leq\left(8+2\pi^{1/2}\right)2^{-J}\TV(f).
\end{align}
By \cref{LinfinityError}, \cref{L1Error}, and the interpolation of $L_p$ spaces, we therefore have
\begin{align*}
\left\|f_J-f\right\|_{L_2(\Om)}&\leq\left\|f_J-f\right\|_{L_1(\Om)}^{1/2}\left\|f_J-f\right\|_{L_{\infty}(\Om)}^{1/2}\leq\left(16+4\pi^{1/2}\right)^{1/2}2^{-J/2}\TV(f)^{1/2}\left\|f\right\|_{L_{\infty}(\Om)}^{1/2}.
\end{align*}
This completes the proof.

\subsection{Proof of \cref{Th4}}\label{PfTh4}

The proof of \cref{Th4} uses the following lemma on the Bessel property of $\left\{\phi(2^J\cdot-\bk):\bk\in\bOm\right\}$.

\begin{lemma}\label{Lemma2} Let $\phi$ be the tensor product piecewise linear B-spline. For each $J\in\N$, we have the followings.
\begin{enumerate}
\item For $u\in L_2(\Om)$, we have
\begin{align}\label{Bessel}
\sum_{\bk\in\bOm}\left|\left\la u,\phi(2^J\cdot-\bk)\right\ra\right|^2\leq\f{4}{2^{2J}}\left\|u\right\|_{L_2(\Om)}^2.
\end{align}
\item For $\bu\in\ell_2(\bOm)$, we have
\begin{align}\label{AdjointBessel}
\left\|\sum_{\bk\in\bOm}\bu[\bk]\phi(2^J\cdot-\bk)\right\|_{L_2(\Om)}^2\leq\f{4}{2^{2J}}\left\|\bu\right\|_{\ell_2(\bOm)}^2.
\end{align}
\end{enumerate}
\end{lemma}

\begin{pf} Note that $0\leq\phi\leq1$, and the direct computation shows that $\left\|\phi\right\|_{L_2(\R^2)}^2\leq1$ and $\left\|\phi(2^J\cdot-\bk)\right\|_{L_2(\R^2)}^2\leq2^{-2J}$. Then by the Schwartz inequality, we have
\begin{align*}
\left|\left\la u,\phi(2^J\cdot-\bk)\right\ra\right|\leq2^{-J}\left(\int_{\Om}\left|u(\x)\right|^21_{\bS_{\bk}}(\x)\rd\x\right)^{1/2}
\end{align*}
where $\bS_{\bk}=\prod_{j=1}^2[2^{-J}(k_j-1),2^{-J}(k_j+1)$. Since $\sum 1_{\bS_{\bk}}=4$, we have
\begin{align*}
\sum_{\bk\in\bOm}\left|\left\la u,\phi(2^J\cdot-\bk)\right\ra\right|^2\leq2^{-2J}\int_{\Om}\left|u(\x)\right|^2\left(\sum_{\bk\in\bOm}1_{\bS_{\bk}}(\x)\right)\rd\x=\f{4}{2^{2J}}\left\|u\right\|_{L_2(\Om)}^2,
\end{align*}
which proves \cref{Bessel}.

For \cref{AdjointBessel}, let $v\in L_2(\Om)$. We have
\begin{align*}
\left\la\sum_{\bk\in\bOm}\bu[\bk]\phi(2^J\cdot-\bk),v\right\ra=\sum_{\bk\in\bOm}\bu[\bk]\left\la\phi(2^J\cdot-\bk),v\right\ra.
\end{align*}
By the Schwartz inequality, we have
\begin{align*}
\sum_{\bk\in\bOm}\left|\bu[\bk]\right|\left|\left\la\phi(2^J\cdot-\bk),v\right\ra\right|\leq\left(\sum_{\bk\in\bOm}\left|\bu[\bk]\right|^2\right)^{1/2}\left(\sum_{\bk\in\bOm}\left|\left\la\phi(2^J\cdot-\bk),v\right\ra\right|^2\right)^{1/2}.
\end{align*}
By \cref{Bessel}, we further have
\begin{align*}
\left|\left\la\sum_{\bk\in\bOm}\bu[\bk]\phi(2^J\cdot-\bk),v\right\ra\right|\leq\left(\sum_{\bk\in\bOm}\left|\bu[\bk]\right|^2\right)^{1/2}\f{2}{2^J}\left\|v\right\|_{L_2(\Om)}.
\end{align*}
Since $v\in L_2(\Om)$ is arbitrary, we have \cref{AdjointBessel} by the converse of H\"older's inequality with $p=q=2$ (e.g. \cite{Folland1999}). This completes the proof.\qquad$\square$
\end{pf}

\begin{poth4} As in \cref{Th5}, it suffices to prove \cref{Th4} for $f\in C^{\infty}(\Om)$ with $f(\x)\in[0,M]$ for $\x\in\Om$, by the standard density argument. Note that we have $d=2$ and $\left|\bOm\right|=2^{2J}$. In addition, since $f(\x)\in[0,M]$ for $\x\in\Om$, $\bsf[\bk]\in[0,M]$ for $\bk\in\bOm$. What is left is to determine $C_{\bsf}$ and $b\in[0,1)$ such that $\left\|\Na\bsf\right\|_1\leq C_{\bsf}2^{2bJ}$. From \cref{Sampling1}, we have
\begin{align*}
\bsf[\bk+\e_j]-\bsf[\bk]=2^J\left\la f,\varphi_{J,\bk+\e_j}\right\ra-2^J\left\la f,\varphi_{J,\bk}\right\ra=2^J\left(\left\la f,\varphi_{J,\bk+\e_j}-\varphi_{J,\bk}\right\ra\right).
\end{align*}
Then, since we have
\begin{align*}
1_{[0,1)}(\cdot-1)-1_{[0,1)}=-\f{\rd}{\rd x}B_2(\cdot-1),
\end{align*}
the direct computations gives
\begin{align*}
\bsf[\bk+\e_j]-\bsf[\bk]=\left\la\f{\p f}{\p x_j},\wt{\varphi}_{j,J,\bk+\e_j}\right\ra,
\end{align*}
where $\wt{\varphi}_j$ is defined as
\begin{align*}
\wt{\varphi}_1(\x)=B_2(x_1)1_{[0,1)}(x_2)~~~\text{and}~~~\wt{\varphi}_2(\x)=1_{[0,1)}(x_1)B_2(x_2).
\end{align*}
Together with the fact that $\wt{\varphi}_j(\cdot+\bk)$ forms a partition of unity for each $j=1,2$, we have
\begin{align*}
\left\|\Na\bsf\right\|_1\leq2^J\int_{\Om}\rd\left\|\na f\right\|_1=2^J\TV(f).
\end{align*}
By setting $C_{\bsf}=\TV(f)$ and $b=1/2$, we have
\begin{align*}
\f{1}{2^{2J}}\left\|\bu^{\La}-\bsf\right\|_{\ell_2(\bOm)}^2\leq\f{64}{3}M^2\left(4+3\sqrt{5(4a+1)\TV(f)\left(1+\TV(f)\right)}\right)\rho^{-1/2}2^{-J/2}\sqrt{8J^3}+\f{16}{3}\eta^2\\
\end{align*}
for some $a\geq1/2$. Hence, we establish \cref{MainResult:Recall} by setting
\begin{align*}
\wt{c}=\f{128}{3}M^2\left(4+3\sqrt{5(4a+1)\TV(f)\left(1+\TV(f)\right)}\right)\sqrt{2},
\end{align*}
with probability at least $1-2^{-2J}$.

For \cref{Coro1:Result}, notice that
\begin{align*}
\left\|u_J^{\La}-f\right\|_{L_2(\Om)}\leq\left\|u_J^{\La}-f_J\right\|_{L_2(\Om)}+\left\|f_J-f\right\|_{L_2(\Om)}.
\end{align*}
More precisely,
\begin{align*}
\left\|u_J^{\La}-f\right\|_{L_2(\Om)}^2\leq2\left(\left\|u_J^{\La}-f_J\right\|_{L_2(\Om)}^2+\left\|f_J-f\right\|_{L_2(\Om)}^2\right).
\end{align*}
By \cref{AdjointBessel} in \cref{Lemma2}, we have
\begin{align*}
\left\|u_J^{\La}-f_J\right\|_{L_2(\Om)}^2\leq\f{4}{2^{2J}}\left\|\bu^{\La}-\bsf\right\|_{\ell_2(\bOm)}^2
\end{align*}
In addition, by \cref{InterpolationApproximation} in \cref{Th5}, we have
\begin{align*}
\left\|f_J-f\right\|_{L_2(\Om)}^2\leq\left(16+4\pi^{1/2}\right)\TV(f)M2^{-J}.
\end{align*}
Hence, we obtain \cref{Coro1:Result} with probability at least $1-2^{-2J}$ by setting
\begin{align*}
C_1=8\wt{c},~~C_2=\f{128}{3},~~\text{and}~~C_3=\left(32+8\pi^{1/2}\right)\TV(f)M,
\end{align*}
all of which are independent of $J$, $\rho$, or $\eta$. This completes the proof.\qquad$\square$
\end{poth4}

\section{Conclusion}\label{Conclusion}

In this paper, we establish an approximation property of total variation minimization from incomplete data. Our error analysis is based on the combination of the uniform law of large numbers and the estimation for its involved covering number of a hypothesis space of the solution. Finally, we further connect our error analysis to the approximation of data with a sparse gradient and the approximation of underlying two dimensional $BV$ functions. For the future work, we plan to establish an approximation from the data on the graph via a graph total variation (e.g. \cite{D.I.Shuman2013}). We may also consider the approximation analysis of the nonlocal total variation (e.g. \cite{X.Zhang2010a}) for the missing data restoration.

\section*{References}

\bibliographystyle{cas-model2-names}

\end{document}